\documentclass[final,onefignum,onetabnum]{siuro210301}

\usepackage{lipsum}
\usepackage{amsfonts}
\usepackage{graphicx, graphics, adjustbox}
\usepackage{epstopdf}
\usepackage{hyperref}
\usepackage{algorithmic}
\usepackage{multicol}
\usepackage{todonotes}
\usepackage{enumitem}
\usepackage{makecell}
\PassOptionsToPackage{table}{xcolor}
\usepackage{tabularray}

\ifpdf
  \DeclareGraphicsExtensions{.eps,.pdf,.png,.jpg}
\else
  \DeclareGraphicsExtensions{.eps}
\fi

\usepackage{enumitem}
\setlist[enumerate]{leftmargin=.5in}
\setlist[itemize]{leftmargin=.5in}

\newsiamremark{remark}{Remark}
\newsiamremark{hypothesis}{Hypothesis}
\crefname{hypothesis}{Hypothesis}{Hypotheses}
\newsiamthm{claim}{Claim}

\headers{A Framework for Approximating Perturbed Optimal Control Problems}{ R. Link and E.Ebbighausen}

\title{A Framework for Approximating Perturbed Optimal Control Problems}

\author{Riley Link \thanks{Creighton University
  (\email{rileylink03@gmail.com}).} \and
  Ethan Ebbighausen \thanks{UNC-CH
  (\email{etebbighausen@gmail.com}).}
}

\dedication{\small\textit{Project advisor: A.Alexanderian\thanks{North Carolina State University (\email{aalexan3@ncsu.edu}).}, J.Hart\thanks{Sandia National Laboratories (\email{joshart@sandia.gov}).}} }

\usepackage{amsopn}

\ifpdf
\hypersetup{
  pdftitle={A Framework for Approximating Perturbed Optimal Control Problems},
  pdfauthor={R. Link and E. Ebbighausen}
}
\fi

\begin{document}

\maketitle

\begin{abstract}
We consider trajectory optimal control problems in which parameter uncertainty limits the applicability of control trajectories computed prior to travel. Hence, efficient trajectory adjustment is needed to ensure successful travel. However, it is often prohibitive or impossible to recalculate the optimal control in-transit due to strict time constraints or limited onboard computing resources. Thus, we propose a framework for quick and accurate trajectory approximations by using post-optimality sensitivity information. This allows the reduction of uncertain parameter space and an instantaneous approximation of the new optimal controller while using sensitivity data computed and stored pretransit.

\end{abstract}

\section{Introduction} \label{sec:intro}

Optimal control problems pose physically difficult questions as optimization problems. Given a physical system, often constrained by differential equations, a control function is found such that an objective or cost function is minimized. This method is commonly used in aerospace applications such as the famous example of determining when to apply thrusters on the Apollo 11 mission \cite{Apollo11}. 

In aerospace applications, the governing equations typically contain parameters that are subject to uncertainty. Examples include environmental parameters such as atmospheric density or physical parameters such as drag and lift coefficients. This parameter uncertainty restricts the usefulness of the precomputed solution to the optimal control problem. A small change in lift or drag forces due to, for example, heat stress causing deformation of a shuttle, can greatly influence how an aircraft maintains trajectory. For small changes, feedback controllers may maintain trajectory \cite{feedbackZhang}. In more extreme cases, an optimal trajectory must be recomputed. We consider these in-transit computations. 

Recalculation of the optimal path is computationally expensive, time-consuming, and lacks any guarantee of obtaining a favorable solution. If a new optimal path is required to avoid failure, a quick and efficient approximation of the optimal path is ideal \cite{HuangRendezvous,ESA}. 

One approach to this computation is presented in \cite{f18Control}. Their adaptive flight control scheme for an F-18 aircraft focuses on in-flight failures that they simulate, and it helps dampen the affects of parameter uncertainties. We present a new method to approach this approximation. 

Section \ref{sec: method} outlines our methods in generality, allowing for its application to a wide range of problems. Our first goal is to reduce the number of significant parameters using global sensitivity analysis. Then, rather than totally recomputing a solution, we focus on approximating a solution, generalizing a method presented in \cite{alexanderian2023new}. By storing and interpolating sensitivity information, which requires less memory due to parameter reduction, we may further reduce in-transit computation time without significant sacrifices of accuracy. These three steps, (i) parameter reduction, (ii) sensitivity information interpolation, and (iii) approximating using a general method from \cite{alexanderian2023new}, individually constitute ways to improve computation time given uncertainty, and we consider the efficacy of each. Hence, our approach is to adapt the optimal trajectory dynamically based on sensitivity information. Typical feedback control which relies on reactive adjustments, works well for small deviations but may struggle with large deviations, and so this new approach can handle those larger perturbations given its adjustments made using precomputed sensitivity information.

Thus, we want to emphasize that our method is intended to complement, rather than replace, existing feedback control techniques, which remain well-suited for handling small deviations. Additionally, this approach represents a singular strategy within optimal control, and its applicability should be evaluated based on the specific problem at hand. Users should consider its usefulness in conjunction with other established methods to determine the most effective solution for their particular application.

Section \ref{sec:examples} applies this solution method to a model Space Shuttle Trajectory problem from \cite{Betts}. There we solve for the optimal control angle which maximizes the longitude traveled in a model space shuttle's reentry. This provides an in-depth example that demonstrates the strengths and limitations of each part of the procedure.

\section{Generalized Framework for Approximation of the Optimal Control}
\label{sec: method}

As noted above, our framework breaks down into 3 main steps.
\begin{enumerate}[label=(\roman*)]
        \item Conduct global sensitivity analysis to understand which parameters we may remove from consideration.
        \item Interpolate post-optimality local sensitivities to store onboard.
        \item Use the interpolated data to make quick and accurate controller approximations for perturbed parameters.
\end{enumerate}

In particular, we perturb important parameters and solve this new optimal control problem by relying on a method of using the sensitivity information to construct a differential equation whose solution is the new optimal control. 

\subsection{Problem Formulation} \label{sec: Problem}

We consider a general trajectory control governed by a system of ordinary differential equations of the form,
\begin{equation}
    \begin{cases} 
      \boldsymbol{\dot{x}}(t)=\boldsymbol{f}(t,\boldsymbol{x},u;\boldsymbol{p}) & t\in(0,T),\\
      \boldsymbol{x}(0)=\boldsymbol{x}_0,\\
\end{cases}
\end{equation}
where $t$ is time, with final time $T>0$, $\boldsymbol{x}:[0,T]\rightarrow\mathbb{R}^n$ is the state vector with initial condition $\boldsymbol{x}_0$, $u:[0,T]\rightarrow\mathbb{R}$ is the trajectory controller, and $\boldsymbol{p}$ is the parameter vector containing $P$ parameters. To normalize parameter changes, we non-dimensionalize parameters to
\begin{equation} \label{eq:ND parameter}
p_i = (1 + \beta_0 \theta_i) \bar{p}_i, \quad \text{for } i = 1, 2, \ldots, P
\end{equation}

where $\bar{p_i}$ is the nominal value of the parameter, and $\theta_i$ is the reparameterization. The factor $\beta_0$ determines the amount of change taken in these parameters. For example, if $\beta_0=0.1$ this simulates a maximum change of $10\%$ from nominal parameter values.

To find the optimal control $u^*(t)$, we minimize a cost function,
\begin{equation}
\displaystyle\min_{u} J(\boldsymbol{x}(t),u(t),t;\boldsymbol{\theta}).
\end{equation}
and impose inequality constraints and terminal state equality constraints.

% a set of $r$ inequality constraints, $\boldsymbol{\eta}\in\mathbb{R}^r$, and a set of $s$ terminal state equality constraints, $\boldsymbol{\Gamma}\in\mathbb{R}^s$,
% \begin{equation}
%     \begin{cases}
%         \boldsymbol{\eta}(t,\boldsymbol{x},u;\boldsymbol{\theta}) \leq 0 & t\in (0,T) \\
%         \boldsymbol{\Gamma}(\boldsymbol{x}(T))=0.
%     \end{cases}
% \end{equation}

Some optimal control solvers strictly enforce these constraints through spectral methods, but constraints may also be imposed through penalty terms in the cost function $J$. The latter method is less accurate, but is useful in situations with low computation power such as the motivating in-transit case. 

To solve this problem computationally, we discretize the time interval $[0,T]$ as $t_0,t_1,...,t_{N}$ for a chosen $N\in\mathbb{N}$. We also linearly discretize the controller by storing the function $\boldsymbol{u}_N(t)$ as a vector $\left[ u_0,u_1,\ldots,u_{N}\right]$ where 
\begin{equation} \label{eq:discretization}
   \boldsymbol{u}_N(t) = \sum_{i=0}^{N} u_i \phi_{i}(t)
\end{equation}
for basis functions $\phi_{i}(t)$ given by
\begin{equation}
    \phi_{i}(t) = \begin{cases} 
     \frac{t-t_{i-1}}{t_{i}-t_{i-1}}& t \in [t_{i-1},t_{i}] \\
     \frac{t_{i+1}-t}{t_{i+1}-t_{i}}& t \in [t_{i},t_{i+1}] \\
     0 & \text{otherwise}.
    \end{cases}  
\end{equation}
% This discretization also provides local magnitudes $u_i$ for sensitivity information as discussed in subsection \ref{sec:HDSA}.

\subsection{Dimension Reduction}

The various uncertain parameters may influence the problem at highly different degrees. Sensitivity information attempts to capture this information. In the present work, this sensitivity analysis helps us locate those parameters most important to the problem and thus we are locating parameters which have less influence on the optimal control. Fixing unimportant parameters at their nominal values greatly reduces the complexity of the problem when considering possible parameter changes and uncertainty estimation.

We consider the parameters $(\theta_1,...,\theta_{N_{p}}) \in \Theta = \Theta_1 \times \Theta_2 \times \cdots \times \Theta_{P} \subset \mathbb{R}^{P}$, for a space of parameters $\Theta$ on which we define the probability law 

\vspace{-4mm}
\begin{equation}
\mu(d\boldsymbol{\theta}) = \displaystyle \prod^{N_p}_{j=1} \pi_j(\theta_j)d\theta_j
\end{equation} 
for the parameters. We wish to compute Sobol indices, which measure the amount of output variance that is caused by variance in input parameters \cite{SOBOL2001}. It is important to note that Sobol indices come in three different forms, main, interaction, and total. The main Sobol index is simply the ratio of the variance of a single input, say the parameter $\theta_j$, by the total variance of the output, in our case the optimal controller $\boldsymbol{u}$. Thus, this main index measures the direct effects of the input variable and does not account for any interaction effects with other variables. These interaction indices capture how the combined variations of input parameters like $\theta_i$ and $\theta_j$ affect the model output. The total Sobol index captures the contribution of an input variable to the output variance, accounting for both its main effects and all interaction effects involving the input variable.

These total Sobol indices provide a clear ranking of how important input parameters are. Total Sobol indices range from 0 to 1, yet will not sum to 1 given the overlapping of interaction terms. These indices would allow the most unimportant parameters to be reduced to their nominal values and ignored for perturbation analysis. However, computing them can be expensive. Therefore, we bound Sobol indices using derivative-based global sensitivity measures (DGSMs) as described in \cite{cleaves2021global}. 

Note that, if Sobol indices could be precisely computed, the main Sobol indices would provide a more valid reduction of parameter space. This is because, when setting a reduced parameter to its nominal value, we are then ignoring all interaction effects, leaving only the direct effects captured by the main Sobol indices. As such, the main Sobol indices would better describe how important each parameter is and whether it should be retained in the model. 

However, since we are working with upper bounds on Sobol indices, the total index also provides a valid measure of whether or not we should retain or discard a parameter. The total index accounts for both direct and interaction effects, providing a comprehensive measure of a parameter's maximum possible contribution to output variance. For example, given a input parameter with a small upper bound on its total index, we can then safely say that even its maximum possible contribution (direct and interaction effects) to the output variance is negligible. 

Nevertheless, the use of upper bounds comes with limitations. While we can confidently identify parameters with negligible total contributions, if given a large upper bound, i.e. an upper bound above 1, we cannot infer the relative importance of the parameter and its direct versus interaction effects. Thus, we leave these parameters with unclear significance to be further analyzed through the perturbation-based methods in our framework.

The structure of computing the upper bounds on total Sobol indices is as follows. Let $\boldsymbol{u}: \Theta \rightarrow \mathbb{R}^{N+1}$ be a vector-valued function. The derivative-based global sensitivity measure of $\boldsymbol{u}$ with respect to parameter $j$ is defined as

\vspace{-0.5cm}

\begin{equation}\label{eq:DGSMs}
    \mathcal{N}_j(u) =  \sum^{N}_{i=1} \int_{\Theta} \left(\dfrac{\partial u_i}{\partial \theta_j} \right)^{2} \mu (d\boldsymbol{\theta}).
\end{equation}

\medskip

These DGSMs may be use to generate an upper bound on the Sobol indices. In the simple case of a uniform random variable, these connect to the Sobol indices via the following theorem. More complex and general versions may be found in \cite{cleaves2021global}.

\begin{theorem}
    Let $\boldsymbol{u}$ be a vector-valued random variable and let $\Gamma$ be its covariance matrix. Assume the inputs $\theta_i$ are independent and identically distributed uniform random variables $\mathcal{U}(a,b)$. Then, 
    \begin{equation}\label{eq:sobolupperbound}
        \mathbb{S}_{j}^{tot}(\boldsymbol{u}) \leq c_{j}^{p} \frac{\mathcal{N}_{j}(\boldsymbol{u})}{Tr(\boldsymbol{\Gamma})}
    \end{equation}
    where $c_{j}^{p} = \frac{(b-a)^{2}}{\pi^2}$.
\end{theorem}

Using this inequality to bound the Sobol indices from above, the DGSMs can show when parameters are unimportant and may be ignored for small perturbations or totally removed from perturbation analysis. This reduces the computational expense of solving updated problems. 

Numerically, the DGSMs may be computed via a Monte Carlo method and the terms $\frac{\partial u_i}{\partial \theta_j}$ may computed via Hyper-Differential Sensitivity Analysis as seen in Section \ref{sec:HDSA}. In particular, a quasi-Monte Carlo (QMC) method employing low-discrepancy sequences to sample points, as presented in \cite{KUCHERENKO20091135}, is shown to converge faster than a standard Monte Carlo method when computing DGSMs. It is also possible to vectorize the process by interchanging the QMC summation and DGSM summation. With this interchange we have, 

\begin{equation} \label{eq:numericalDGSM}
\mathcal{N}_j(\boldsymbol{u}) \approx \frac{1}{M}  \sum^M_{k=1} \sum^N_{i=1} \left(\dfrac{\partial u^*_i}{\partial \theta_{j,k}} \right)^2.
\end{equation}

Using equation (\ref{eq:numericalDGSM}), upper bounds on Sobol indices for each parameter $\theta_i$ can be calculated. By ignoring the unimportant parameters the complexity of the optimal control problem is reduced and when coupled with Hyper-Differential Sensitivity Analysis, the computational cost for this reduction is feasible.

\subsubsection{Hyper-Differential Sensitivity Analysis} \label{sec:HDSA}

Hyper-Differential Sensitivity Analysis (HDSA) is a local sensitivity method that considers the change in the solution of an optimal control problem with respect to perturbations of parameters appearing in the model \cite{hart2023sensitivity}. It encapsulates this change by computing the derivatives $\frac{\partial \boldsymbol{u}^*}{\partial \theta_i}$, where $\boldsymbol{u}^*$ is the discretized optimal solution and $\theta_i$ are relevant parameters. It is possible to use numerical differentiation to compute these, but it diminishes the error to use implicit differentiation as in \cite{hart2023sensitivity}. 

Considering $\nabla_u J(u^*(\theta),\theta)=0$, implicit differentiation with respect to $\theta$ gives 
\begin{equation}
    \left(\frac{\partial^{2} J}{\partial u_i \partial u_j}\right) \left(\frac{\partial u^{*}_{j}}{\partial \theta_k}\right) + \left(\frac{\partial^{2} J}{\partial u_i \partial \theta_k}\right) = 0, \quad \text{for } i, j = 0, \ldots, N \text{ and } k = 1, \ldots, P
\end{equation}

\vspace{-1mm}
or 
\vspace{-5mm}

\begin{align} \label{eq:D=-HB}
    &\boldsymbol{D}=-\boldsymbol{H}^{-1}\boldsymbol{B}
\end{align}
where $\boldsymbol{H}=(\frac{\partial^{2} J}{\partial u_i \partial u_j})_{ij}$ is the Hessian of $J$ with respect to $u$, which is assumed to be invertible at the optimal solution, $\boldsymbol{B}=(\frac{\partial^{2} J}{\partial u_i \partial \theta_j})_{ij}$ is a matrix of mixed partials with respect to the controller $\boldsymbol{u}$ and parameters $\boldsymbol{\theta}$, and $\boldsymbol{D}$ is the matrix of sensitivities $\frac{\partial \boldsymbol{u}^{*}_{i}}{\partial \theta_k}$. Since $\boldsymbol{u}$ is discretized, the $\boldsymbol{D}$ matrix consists of local sensitivities of the discretized $\boldsymbol{u}$ with respect to each parameter $\boldsymbol{\theta}$. The entries of $\boldsymbol{D}$ are used for calculating DGSMs and for approximating $\boldsymbol{u}^*$ as in the next section and the Hessian and mixed partials matrix may themselves be computed by standard numerical differentiation techniques such as complex step differentiation or sensitivity equations.

\subsubsection{Sensitivity Equations} \label{sec:senseq}
The calculation of the gradient $\frac{\partial J}{\partial u_i}$ of the cost function $J(\boldsymbol{x}(t),\boldsymbol{u}(t),t;\boldsymbol{\theta})$, is important in both HDSA and the convergence and precision of the optimization process. Thus, employing an exact calculation rather than a numerical approximation is crucial. We choose to achieve this through the use of sensitivity equations.

Given a system of ODEs with state variables $x_i$ and a discretized controller $\boldsymbol{u}$, sensitivity of the state variables with respect to the controller can be calculated via sensitivity equations. This is done by treating the discretized controller $\boldsymbol{u}$ as a parameter in the ODEs such that
\begin{align}
\begin{split}
    \boldsymbol{x}'(t;\boldsymbol{u})&=\boldsymbol{f}(t,\boldsymbol{x}(t;\boldsymbol{u});\boldsymbol{u}) \\
    \boldsymbol{x}(0)&=\boldsymbol{x}_0(\boldsymbol{u}).
\end{split}
\end{align}
By treating the controller this way it is possible to differentiate with respect to $\boldsymbol{u}$
\begin{equation}
    \frac{\partial x_i'}{\partial u_j}=\frac{\partial \boldsymbol{f}}{\partial x_i}\frac{\partial x_i}{\partial u_j}+\frac{\partial \boldsymbol{f}}{\partial u_j}.
\end{equation}
The derivatives $\frac{\partial x_i}{\partial u_j}$ can be used to calculate a more accurate gradient, $\frac{\partial J}{\partial u_j}$, when differentiating the cost function directly. Appendix \ref{appendix:senseq} shows calculations of sensitivity equations for the example problem shown in section \ref{sec:examples}.

\subsection{Approximating the Optimal Control Problem}
In case of parameter change during flight (either an actual change occurs or parameters are remeasured at better accuracy), it is necessary to shift the precomputed trajectory to match the physical system, but the time pressure and computational restrictions of in-flight computation impedes the ability to recompute an optimal control, especially in the case of a highly non-convex objective function. Instead, we shift to approximation of the optimal control. In particular, we focus on a method relying on the sensitivity information described in the previous section. 

\subsubsection{Generalizing an Approximation Method}
From subsection \ref{sec:HDSA}, we obtain the matrix $\boldsymbol{D} = (\frac{ \partial u^{*}_{i}}{\partial \theta_j})_{ij}$. For small changes in $\boldsymbol{\theta}$ or changes with respect to less-important parameters, we may apply a first-order Taylor series approximation,

\vspace{-5mm}
\begin{equation}
\boldsymbol{u}^{*}(\boldsymbol{\theta}_1) = \boldsymbol{u}^{*}(\boldsymbol{\theta}_0) + \boldsymbol{D}(\boldsymbol{u}^{*}(\boldsymbol{\theta}_0),\boldsymbol{\theta}_0)(\boldsymbol{\theta}_1-\boldsymbol{\theta}_0) + \mathcal{O}(||\boldsymbol{\theta}_1 - \boldsymbol{\theta}_0||_{2}),
\end{equation}

for nominal parameters $\boldsymbol{\theta}_0$ and perturbed parameters $\boldsymbol{\theta}_1$.
This approximation is already relatively accurate. Computing the derivative information $\boldsymbol{D}$ offline and storing it in vehicle systems allows for a rapid approximation. Higher order approximations may greatly improve accuracy if derivatives are highly accurate as well.

The inaccuracy of the linear approximation depends on the curvature of the optimal controller $\boldsymbol{u}^{*}$ with respect to $\boldsymbol{\theta}$. This problem may be addressed by stepping over $\boldsymbol{\theta}$ as presented in \cite{alexanderian2023new}. This process first parameterizes the step from $\boldsymbol{\theta}_0$ to $\boldsymbol{\theta}_1$ as $\boldsymbol{\theta}(t) = \boldsymbol{\theta}_0 + t(\boldsymbol{\theta}_1 - \boldsymbol{\theta}_0)$. Then, 

\vspace{-2mm}
\begin{equation}
\begin{aligned}
\frac{\partial \boldsymbol{u}^{*}}{\partial t} &= \frac{\partial \boldsymbol{u}^{*}}{\partial \boldsymbol{\theta}} \frac{\partial \boldsymbol{\theta}}{\partial t} \\
&=\boldsymbol{D}(\boldsymbol{u}^{*}(\boldsymbol{\theta}(t)),\boldsymbol{\theta}(t)) (\boldsymbol{\theta}_1-\boldsymbol{\theta}_0). 
\end{aligned}
\end{equation}

This creates a differential equation for $\boldsymbol{u}$ which may be evaluated to obtain $\boldsymbol{u}^{*}(\boldsymbol{\theta}_1) = \boldsymbol{u}^*(\boldsymbol{\theta}(1))$. In \cite{alexanderian2023new}, this is presented via a time-stepping, or modified forward-Euler method. We pick some $M \in \mathbb{N}$, set $h = \frac{1}{M}$ and $t_{m} = mh$, and apply the recursion 

\vspace{-3mm}
\begin{equation}
\label{eq:FE}
\boldsymbol{u}^{*}_{m+1} = \boldsymbol{u}^{*}_{m} + h\boldsymbol{D}(\boldsymbol{u}^{*}_{m},\boldsymbol{\theta}(t_{m}))(\boldsymbol{\theta}_{1}-\boldsymbol{\theta}_{0})
\end{equation}
so $\boldsymbol{u}^{*}(\boldsymbol{\theta_1}) \approx \boldsymbol{u}^{*}_{M}$. Since this method does require re-computation of the derivative at various points, it trades accuracy for some computational speed. This is well-known to have a local truncation error proportional to the square of the step size, or $\mathcal{O}(||\boldsymbol{\theta_1} - \boldsymbol{\theta_0}||_{2}^{2})$. In the case of in transit changes or large perturbations, this accuracy becomes extremely important. However, due to compounding error in calculating derivatives and steps, the number of step sizes that improve the approximation is limited by the system involved. Any number of other numerical differential equation methods may be applied instead of a forward-Euler method, such as higher order Runge-Kutta methods which may help control when the ODE acts stiff at the cost of computation time.

\subsubsection{Interpolation of the Key Jacobian}
In $n^{th}$-order Runge Kutta, the number of calculations of the derivative is on the scale of $\mathcal{O}\left(nMNP\right)$ for $M$ the number of steps, $N$ the discretization number, and $P$ the number of parameters. Thus, the computation of these derivatives throttles our increase of the order for accuracy. Depending on the curvature of the objective, sacrificing some accuracy of the derivatives to improve the order of the approximation may improve the approximation of the controller.

In particular, we focus on approximating the key Jacobian, $\boldsymbol{D} = (\frac{ \partial u^{*}_{i}}{\partial \theta_j})_{ij}$. We precompute the sensitivity matrix, $\boldsymbol{D}$, on a mesh of values in the hypercube $\Theta$, and store grid-call interpolation functions to approximate the derivative at given points. While this can be taxing in terms of computational power prior to an experiment, it drastically improves computation time in-situ. In some cases it improves accuracy due to the compounding error as the number of ODE steps increases. We will refer to the ODE method using these derivatives for inputs as the Interpolated Step method in the example shown in section \ref{sec:examples}.

Ironically, an increased number of steps tends to outweigh higher order methods in terms of accuracy in practice. While several other ODE methods, such as Runge Kutta 3, Adams-Bashforth multistep, Adams-Moulton multistep, and Gaussian Quadrature, were considered, we will focus on the forward Euler method moving forward and show the effectiveness of the interpolated step in section \ref{sec:examples}. Additionally, while it is not implemented in this paper, storing higher-order derivatives to approximate the Jacobian $\frac{\partial u_i}{\partial \theta_j}$ could similarly improve computational efficiency without sacrificing much accuracy.

\section{Numerical Results} \label{sec:examples}
To demonstrate the effectiveness and efficacy of our method, we provide an example of a 2 degree-of-freedom space shuttle re-entry problem. 

Less complex problems, such as the classical Zermelo's navigation problem, were considered but omitted as the small parameter space and ease of solution did not warrant parameter reduction or approximation, and did not display noticeable change with their application. However, with more complex systems, such as the 2 DOF space shuttle re-entry, in-transit parameter shifts can cause catastrophic results and thus require approximations to the optimal control to be made. Furthermore, these problems contain many physical parameters and thus demonstrate the effectiveness of our parameter reduction method and approximation technique.

To achieve our proposed method, we implement this numerically in $\textsc{matlab}$. Algorithm \ref{algorithm} provides a summary of the proposed approach. 

\vspace{4mm}

\begin{algorithm}[H]
\caption{Computation of Mid-Trajectory Approximations}
\label{algorithm}
\begin{algorithmic}[1]
\STATE{Solve the nominal trajectory problem $\left(\boldsymbol{\theta}_0=\boldsymbol{0}\right)$ for states $\boldsymbol{x}^*$ and optimal controller $\boldsymbol{u}^*$}
\STATE{Subdivide the hypercube $\Theta$ into $m$ equally distributed samples. At each sample compute the optimal trajectory and take local hyper-differential sensitivities, $\boldsymbol{D}$ (\ref{eq:D=-HB}), about the new trajectory. These sensitivities form the basis used for interpolation }
\STATE{Compute DGSMs (\ref{eq:DGSMs}) and reduce the parameter space using the Sobol Index upper bounds (\ref{eq:sobolupperbound})}
\STATE{Organize the sensitivities computed across all of $\Theta$ into a key Jacobian of sensitivities for only the important parameters.}
\STATE{Assume a random parameter change, $\boldsymbol{\theta}_1$, at $t_0$ such as $\theta_j\in \mathcal{U}([-1,1])$}
\STATE{Using (\ref{eq:FE}), step across $\boldsymbol{\theta}_0\to\boldsymbol{\theta}_1$ and interpolate the key Jacobian to calculate the sensitivities at each step in order to return the approximate trajectory $\boldsymbol{u}_{IS}$}
\end{algorithmic}
\end{algorithm}

\subsection{Space Shuttle Problem} \label{sec:spaceshuttle}
We consider a 2 degree-of-freedom space shuttle re-entry problem posed in \cite{Betts}. With 7 parameters to consider, this gives a sufficient amount of complexity to examine the method described above. In particular, we note that the parameter reductions and Jacobian interpolation do not substantially affect the accuracy of the approximated optimal control (measured through norm difference of the controls and difference of the associated costs). Further, we display the ability of this approximation method to mimic the new optimal solution while noting some instability in this approximation method.

\subsubsection{Problem Formulation}
The motion of the space shuttle is defined by the following set of equations:

\vspace{3mm}

\begin{minipage}{0.6\textwidth}
\begin{equation}
\boldsymbol{\dot{x}}(t)=
\begin{cases}
    \dot{h}(t)&=v(t)\sin{(\gamma(t))} \\
    \dot{\phi}(t)&= \frac{v(t)}{r}\cos{(\gamma(t))} \\
    \dot{v}(t)&=-\frac{D}{m} - g\sin{(\gamma(t))} \\
    \dot{\gamma}(t)&=\frac{L}{mv(t)} + \cos{(\gamma(t))}\left( \frac{v(t)}{r}-\frac{g}{v(t)} \right)
\end{cases}
\end{equation}
\end{minipage}
\begin{minipage}{0.35\textwidth}
\begin{table}[H]
\begin{center}
    \begin{tabular}{c c}
    \hline
    $h$ & altitude (ft) \\
    $\phi$ & longitude (deg) \\
    $v$ & velocity (ft/sec) \\ 
    $\gamma$ & flight path angle (rad) \\
    \hline
    \end{tabular}
\end{center}
\end{table}
\end{minipage}

\vspace{5mm}

\noindent The initial conditions, terminal conditions, and constraints of the system are as follows:

% \begin{equation} \label{eq:SS initial}
% \boldsymbol{x}(0)=
%     \begin{cases}
%     h(0)&=80000\\
%     v(0)&=2500\\
%     \gamma(0)&= -5\frac{\pi}{180}.
%     \end{cases}
% \end{equation}

% \begin{equation} \label{eq:SS final}
% \boldsymbol{x}(T)=
%     \begin{cases}
%         h(T)&=80000\\ v(T)&=2500\\
%         \gamma(T)&= -5\frac{\pi}{180}.
%     \end{cases}
% \end{equation}

% \begin{equation} \label{eq:SS constraints}
% \boldsymbol{\eta}(t,\boldsymbol{x},u;\boldsymbol{\theta})=
%     \begin{cases}
%         h(t) &\geq 0 \\
%         v(t)&\geq 1 \\
%         |\gamma (t)|&\leq 89\frac{\pi}{180} \\
%         |u(t)| &\leq\frac{\pi}{2}.
%     \end{cases}
% \end{equation}

\begin{minipage}{0.45\textwidth}
\begin{equation} \label{eq:SS initial}
\boldsymbol{x}(0)=
    \begin{cases}
    h(0) &= 260000 \\
    v(0) &= 25600 \\
    \gamma(0) &= -\frac{\pi}{180}
    \end{cases}
\end{equation}
\end{minipage}%
\hfill
\begin{minipage}{0.45\textwidth}
\begin{equation} \label{eq:SS final}
\boldsymbol{x}(T)=
    \begin{cases}
        h(T) &= 80000 \\
        v(T) &= 2500 \\
        \gamma(T) &= -5\frac{\pi}{180}
    \end{cases}
\end{equation}
\end{minipage}

\begin{equation} \label{eq:SS constraints}
\boldsymbol{\eta}(t,\boldsymbol{x},u;\boldsymbol{\theta})=
    \begin{cases}
        h(t) &\geq 0 \\
        v(t) &\geq 1 \\
        |\gamma(t)| &\leq 89\frac{\pi}{180} \\
        |u(t)| &\leq \frac{\pi}{2}
    \end{cases}
\end{equation}

\noindent Appendix \ref{appendix:SSvariables} contains the specific aerodynamic and atmospheric forces on the shuttle.

The goal of the space shuttle re-entry is to maximize longitudinal distance, $\phi(t):=x_2(t)$, across $t\in (0,T)$. This is achieved by finding the optimal angle of attack $u:[0,T]\to \mathbb{R}$ such that the cost function $J(\boldsymbol{x}(t),u(t),t;\boldsymbol{\theta})$ is minimized. Thus let the optimal control problem be formulated as,
\begin{equation} \label{eq:ssMIN}
   \displaystyle\min_{u} J(\boldsymbol{x}(t),u(t),t;\boldsymbol{\theta})=\displaystyle\min_{u} \left[ -x_2(T)+S(\boldsymbol{x}(t))\right]
\end{equation}
where $S(\boldsymbol{x}(t))$ is a function of penalty terms enforcing terminal conditions (\ref{eq:SS final}) and state constraints (\ref{eq:SS constraints}). The full cost function is shown in appendix \ref{appendix:cost}. In particular, we take a weighted penalty for each state and exponentiate to ensure the cost is twice-differentiable.

We analyze seven parameters of interest that are subject to change during flight or to uncertainty in measurement:
\begin{align}
    \boldsymbol{p} = [m,\rho_0,a_0,a_1,b_0,b_1,b_2]
\end{align}
where $m$ is the mass of the space shuttle, $\rho_0$ the initial atmospheric density, and $a_i$ and $b_i$ are dimensionless parameters found in the coefficients of lift and drag. As shown in subsection \ref{sec: Problem}, these parameters are nondimensionalized to $\boldsymbol{\theta}\in\Theta$ according to equation (\ref{eq:ND parameter}).

To solve the problem computationally, we linearly discretize the controller according to equation (\ref{eq:discretization}). Thus the discretized controller is $\boldsymbol{u}:[0,T]\to\mathbb{R}^{N+1}$. A discretization value of $N=10$ or $N=20$ proved sufficient for computations with appropriate interpolation of the controller when necessary.

Sensitivity equations were used for efficient computation of the gradient, $\frac{\partial J}{\partial u_i}$, and the full computations can be seen in appendix \ref{appendix:senseq}. As mentioned previously, the calculation of the gradient through this method also increases the accuracy of the local sensitivities, $\frac{\partial u_i^*}{\partial \theta_j}$, which are found through HDSA as outlined in subsection \ref{sec:HDSA}. 

%A graph of the local sensitivities is shown in figure \ref{}.

\subsubsection{Parameter Space Reduction}\label{sec:SSparamreduc}
We first assume that all nondimensionalied parameters, $\theta_j$, follow a uniform distribution, $\mathcal{U}([-1,1])$, which allows for a maximum 10\% change in the parameters. This assumption is made due to the absence of further information on the atmospheric space of interest or information on the non-dimensional drag and lift coefficients.

We computed DGSM's and upper bounds on Sobol Indices according to equations (\ref{eq:numericalDGSM}) and (\ref{eq:sobolupperbound}) respectively. For a sample of $M=800$ the results are shown in Table \ref{tab:SS DGSM}.

\begin{table}[H]
\footnotesize
\caption{\label{tab:SS DGSM} DGSM's and upper bounds of Sobol Indices for parameters one through seven.}
\begin{center}
\begin{tblr}{
 colspec={||c c c||}, 
 rowsep=3pt, % Adjust row spacing
 hlines,
 row{6}={gray!30}, % Shades row with a_1
 row{5}={gray!30}, % Shades row with b_0
 row{8}={gray!30}  % Shades row with b_2
}
 Parameter \# & DGSM & Upper Bound \\ 
 $m$ & \( 3.22 \times 10^{-3} \) & \( 2.13 \times 10^{-2} \) \\ 
 $\rho_0$ & \( 3.56 \times 10^{-3} \) & \( 2.35 \times 10^{-2} \) \\
 $a_0$ & \( 9.91 \times 10^{-3} \) & \( 6.54 \times 10^{-2} \) \\
 $a_1$ & \( 2.66 \times 10^{0} \) & \( 1.76 \times 10^{1} \) \\
 $b_0$ & \( 6.56 \times 10^{-1} \) & \( 4.33 \times 10^{0} \) \\
 $b_1$ & \( 5.39 \times 10^{-3} \) & \( 3.56 \times 10^{-2} \) \\
 $b_2$ & \( 7.50 \times 10^{-1} \) & \( 4.95 \times 10^{0} \) \\
\end{tblr}
\end{center}
\end{table}

% \begin{table}[H]
% \footnotesize
% \caption{\label{tab:SS DGSM} DGSM's and upper bounds of Sobol Indices for parameters one through seven.}
% \begin{center}
% \begin{tblr}{
%  colspec={||c c c||}, 
%  rowsep=0pt, % Adjust row spacing
%  hlines,
%  row{6}={gray!30}, % Shades row with a_1
%  row{5}={gray!30}, % Shades row with b_0
%  row{8}={gray!30}  % Shades row with b_2
% }
%  Parameter \# & DGSM & Upper Bound \\ 
%  $m$ & 0.00322 & 0.0213 \\ 
%  $\rho_0$ & 0.00356 & 0.0235 \\
%  $a_0$ & 0.00991 & 0.0654 \\
%  $a_1$ & 2.66 & 17.6 \\
%  $b_0$ & 0.656 & 4.33 \\
%  $b_1$ & 0.00539 & 0.0356 \\
%  $b_2$ & 0.750 & 4.95 \\
% \end{tblr}
% \end{center}
% \end{table}

Table \ref{tab:SS DGSM} shows that the mass $m$, atmospheric density $\rho_0$, $a_0$, and $b_1$ parameters are relatively unimportant. Consequently, we focus on perturbations to parameters $a_1$, $b_0$, and $b_2$. Thus we have reduced the total number of parameters to $N_p=3$ and so $\Theta\subset\mathbb{R}^3$, allowing for quicker computations. This process is shown to not diminish accuracy through results in subsection \ref{sec:SSresults}.

\subsubsection{Space Shuttle Results}\label{sec:SSresults}

\begin{center}
    \begin{figure}[H]
        \includegraphics[width=0.95\textwidth, height = 10cm, trim={3mm 6mm 3mm 6mm},clip]{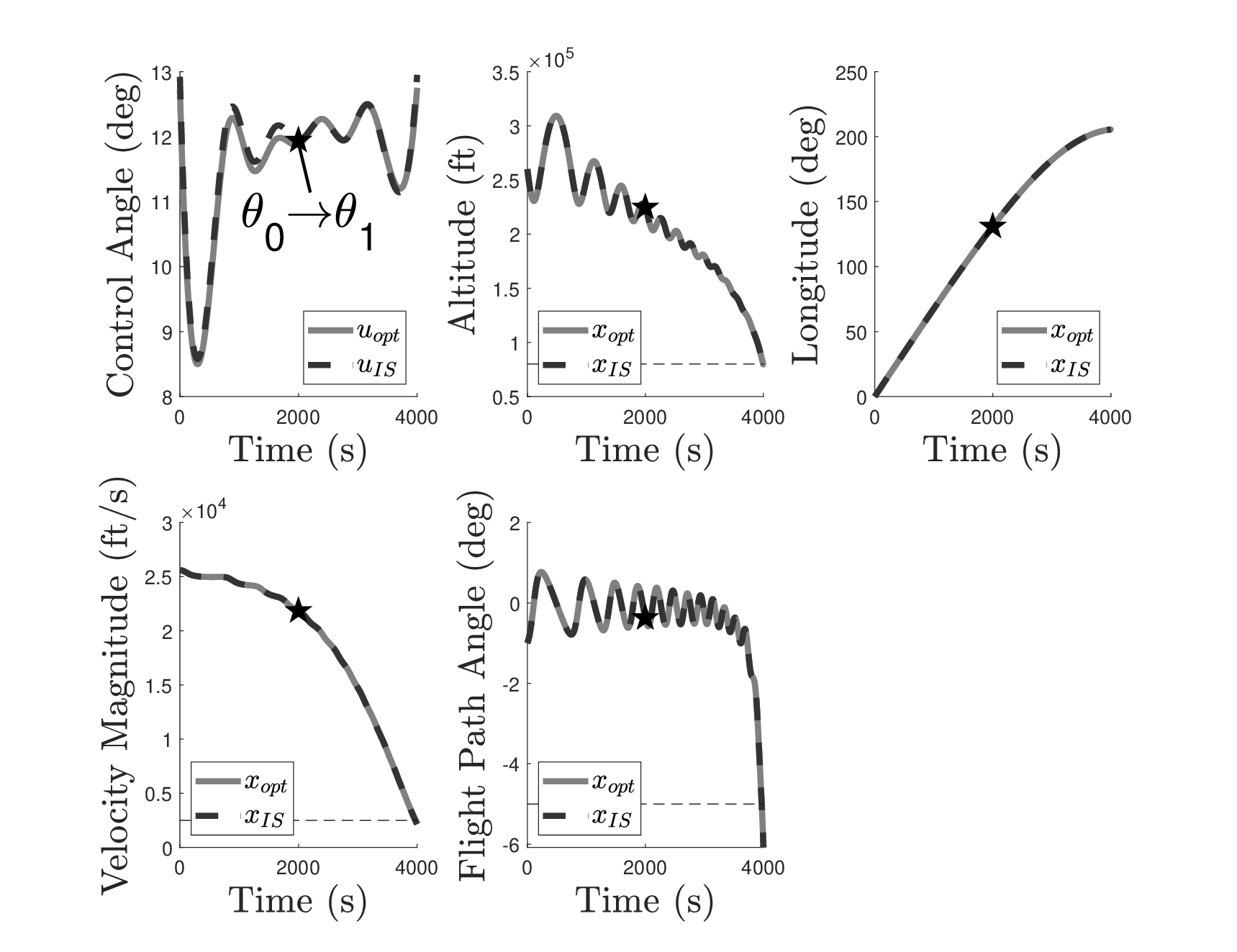}
        \caption{Controller, $\boldsymbol{u}$, and state variables, $\boldsymbol{x}$, of the optimal solution and approximated solution for parameter shift $\boldsymbol{\theta}_0\to\boldsymbol{\theta}_1$ at $t=2000$. The terminal constraints are represented by the thin dashed lines.}
        \label{fig:ssGraph}
    \end{figure}
\end{center}

At nominal parameters $\boldsymbol{\theta}_0=\boldsymbol{0}$, the optimal solution is denoted $u^*$. This is the expected path the space shuttle would follow given no perturbations in parameters. We simulate a parameter change $\boldsymbol{\theta}_0\to\boldsymbol{\theta}_1$ at $t=2000$ in which only the important parameters found in subsection \ref{sec:SSparamreduc} are perturbed. The optimal solution for the parameter change is denoted $u_{opt}$. Figure \ref{fig:ssGraph} shows an approximation using our interpolated step method for a random parameter perturbation $\boldsymbol{\theta}_0\to\boldsymbol{\theta}_1$ at $t=2000$. This simulates a change halfway through the trajectory that requires a new path. In this example the approximation is shown to closely follow the new optimal solution in this parameter shift. However, note that in Figure \ref{fig:ssGraph}, even the optimal solution violates some final state constraints. This is due to the fact that we impose constraints as penalties, allowing for some violations at the expense of an increase in the cost function. 

\vspace{-7mm}

\begin{center}
    \begin{figure}[h]
        \includegraphics[width=0.8\textwidth,height = 6cm, trim={0.5cm 0cm 0.5cm 0.5cm},clip]{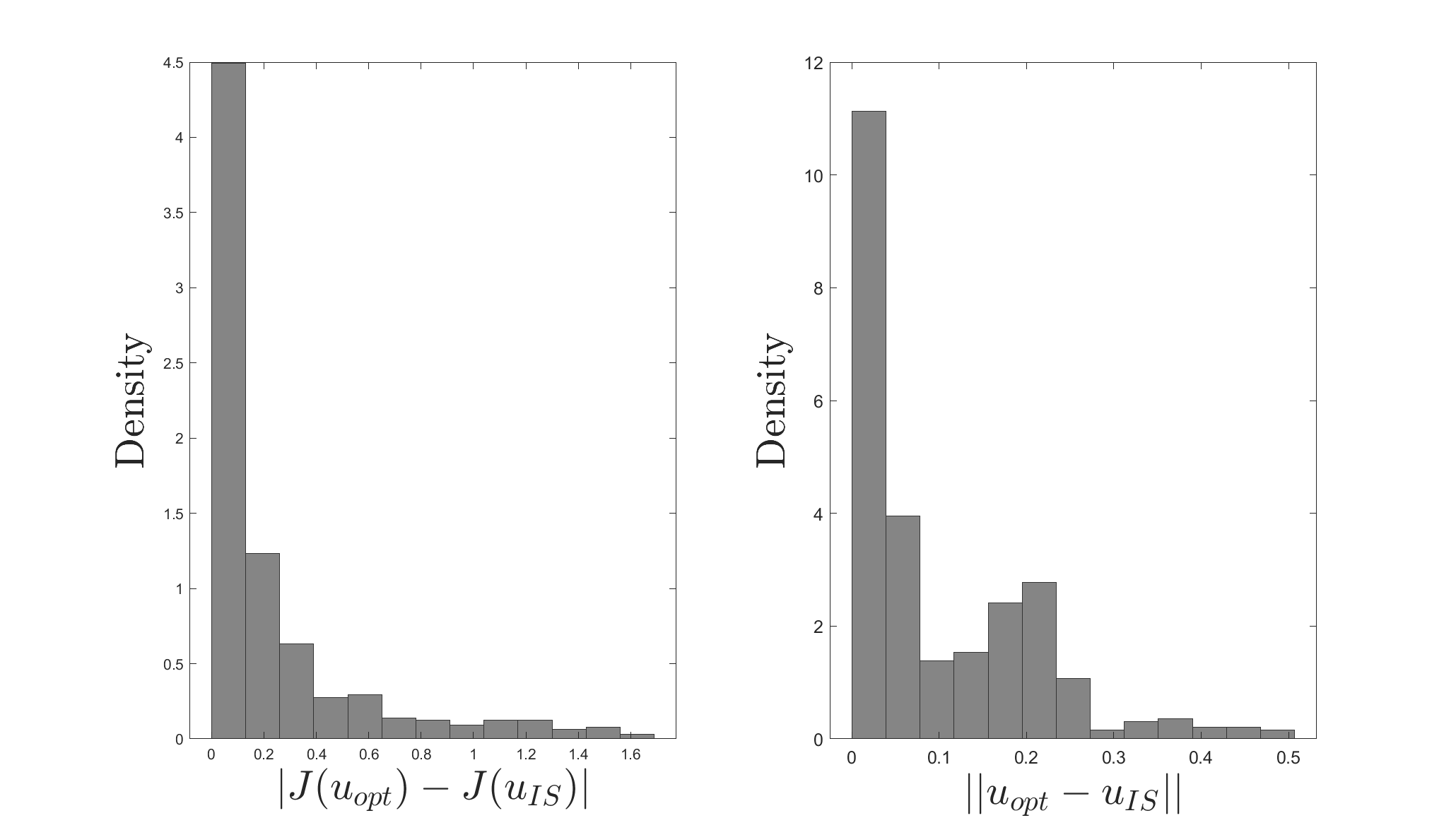}
        \vspace{-0.3cm}
        \caption{Cost difference (left) and error norm (right) between optimal and approximated solutions for 500 parameter perturbations. Note that the cost is on the scale of $2$ and $||u||$ is on the scale of $0.68$. }
        \label{fig:ssHist}
    \end{figure}
\end{center}

After $500$ samples of random parameter perturbations, in which only the important parameters were perturbed, Figure \ref{fig:ssHist} gives histograms of results to display the relative density of approximations. 

\vspace{-4mm}

In Figure \ref{fig:ssHist} we see that both the cost and the norm error skew right, thus implying most interpolated step approximations efficiently minimize cost and accurately approximate the optimal solution for the parameter change. We chose these two error metrics because we noticed that an approximation that minimizes error does not necessarily minimize cost and vice versa. Since we only have control over the costs here, this implies there are multiple paths that lead to similar results, and computationally this implies that this is a highly non-convex problem with many local minima. By using both metrics, a more complete view of the effectiveness of our approximation is shown. 

Based on the scale of the cost and norm we notice that there are outliers for both cost and error measurements. Since we take random parameter perturbations, there are situations when they are beyond the ability of our optimizer (without extra care to find a better initial guess), and the approximation fails in a similar way to the re-optimized result.

Finally, Figure \ref{fig:ssHist} displays a bimodal distribution with a main mode around 0 and a second mode around 0.2. We notice that this secondary mode seems to be a feature of the problem. As more samples are taken, the tail extends and the mode persists. However, we can calculate, 
\begin{align*}
    \mathbb{P}\left(\lVert u_{opt}-u_{IS}\rVert<0.2\right)&=0.84.
\end{align*} 
Thus, there is an 84\% chance that our approximation is below this secondary mode and thus is an excellent error approximation. Additionally, given that we are looking at a considerable 10\% perturbation range, we would expect this probability to only increase with a more local approximation perturbation range of 5\% or less. 

Next, to verify the usefulness of parameter reduction, we use the 500 parameter changes and compare the new optimal solution using all $7$ parameter perturbations, $u_{opt}^{7}$, with the case of only considering the $3$ parameter perturbations identified in subsection \ref{sec:SSparamreduc}, $u_{opt}^{3}$. The following table summarizes the results:
\vspace{-1mm}

\renewcommand{\arraystretch}{1.3} 
\begin{table}[H]
\footnotesize
\caption{\label{tab:opt solutions} Optimal Solutions and Parameter Reduction for 500 Perturbations}
\begin{center}
\begin{tabular}{||c|| c c c||} 
 \hline
 \rule{0pt}{2.7ex}
  & \( J\left(u_{opt}^{7}\right) \) & \( J\left(u_{opt}^{3}\right) \) & \( \left\lVert u_{opt-7p}-u_{opt-3p}\right\rVert \) \\ [0.5ex] 
 \hline
 Mean & \( 1.37 \times 10^0 \) & \( 1.11 \times 10^0 \) & \( 6.21 \times 10^{-2} \) \\ 
 \hline
 Median & \( -1.63 \times 10^0 \) & \( -1.60 \times 10^0 \) & \( 2.11 \times 10^{-2} \) \\
 \hline
\end{tabular}
\end{center}
\end{table}

% \begin{table}[H]
% \footnotesize
% \caption{\label{tab:opt solutions} Optimal Solutions and Parameter Reduction for 500 Perturbations}
% \begin{center}
% \begin{tabular}{||c|| c c c||} 
%  \hline
%  \rule{0pt}{2.7ex}
%   & $J\left(u_{opt}^{7}\right)$ & $J\left(u_{opt}^{3}\right)$ & $\left\lVert u_{opt-7p}-u_{opt-3p}\right\rVert$ \\ [0.5ex] 
%  \hline
%  Mean & 1.3729 & 1.1104 & 0.0621 \\ 
%  \hline
%  Median & -1.6338 & -1.5973 & 0.0211 \\
%  \hline
% \end{tabular}
% \end{center}
% \end{table}

\noindent where $||u||$ is on the scale of $0.68$. The median difference, 0.0211, represents about a $3\%$ change in the controller norm and an even smaller change in the optimal solution's cost. Thus, optimizing with the reduced parameters returns solutions with nearly the exact same cost and almost no difference in path.

Note that in this study, the constraints are encoded as penalties in the objective, so instability represents itself as the discrepancy between the median and mean in the approximation costs, whereby the small norm difference highlights the nearness to the optimizer output. 

Table \ref{tab:LIN cost SS param} considers the linear (or forward-Euler) approximation over the same sample of perturbed parameters sub-scripted as $u_{lin}$, both over 7 parameters and 3 parameters with notation as above. Table \ref{tab:LIN error SS param} compares the norm differences to the reoptimized solution. 
\vspace{-1mm}

\begin{table}[H]
\footnotesize
\caption{\label{tab:LIN cost SS param} 500 Perturbed Parameters Linear Approximation Costs without Reduction v.s. with Reduction}
\begin{center}
\begin{tabular}{||c|| c c c||} 
 \hline
  & \( J\left(u_{opt}^{7}\right) \) & \( J\left(u_{lin}^{7}\right) \) & \( J\left(u_{lin}^{3}\right) \) \\ [0.5ex] 
 \hline
 Mean & \( 1.37 \times 10^0 \) & \( 2.54 \times 10^3 \) & \( 6.83 \times 10^2 \) \\ 
 \hline
 Median & \( -1.63 \times 10^0 \) & \( -9.68 \times 10^{-1} \) & \( -9.42 \times 10^{-1} \) \\
 \hline
\end{tabular}
\end{center}
\end{table}

\vspace{-2mm}

\begin{table}[H]
\footnotesize
\caption{\label{tab:LIN error SS param} 500 Perturbed Parameters Linear Approximation Errors without Reduction v.s. with Reduction}
\begin{center}
\begin{tabular}{||c|| c c c c c||} 
 \hline
  & \( \left\lVert u_{lin}^{7}-u_{lin}^{3}\right\rVert \) 
  & \( \left\lVert u_{opt}^{7}-u_{lin}^{7}\right\rVert \) 
  & \( \left\lVert u_{opt}^{7}-u_{lin}^{3}\right\rVert \) 
  & \( \left\lVert u_{opt}^{3}-u_{lin}^{7}\right\rVert \) 
  & \( \left\lVert u_{opt}^{3}-u_{lin}^{3}\right\rVert \) \\ [0.5ex] 
 \hline
 Mean & \( 4.27 \times 10^{-2} \) & \( 1.48 \times 10^{-1} \) & \( 1.45 \times 10^{-1} \) & \( 1.44 \times 10^{-1} \) & \( 1.37 \times 10^{-1} \) \\ 
 \hline
 Median & \( 4.37 \times 10^{-2} \) & \( 1.21 \times 10^{-1} \) & \( 1.19 \times 10^{-1} \) & \( 1.26 \times 10^{-1} \) & \( 1.17 \times 10^{-1} \) \\
 \hline
\end{tabular}
\end{center}
\end{table}

% \renewcommand{\arraystretch}{1.25}
% \begin{table}[H]
% \footnotesize
% \caption{\label{tab:LIN cost SS param} 500 Perturbed Parameters Linear Approximation Costs without Reduction v.s. with Reduction}
% \begin{center}
% \begin{tabular}{||c|| c c c||} 
%  \hline
%   & $J\left(u_{opt}^{7}\right)$ & $J\left(u_{lin}^{7}\right)$ & $J\left(u_{lin}^{3}\right)$ \\ [0.5ex] 
%  \hline
%  Mean & 1.3729 & $2.54\times 10^3$ & 683.1728 \\ 
%  \hline
%  Median & -1.6338 & -0.9683 & -0.9419 \\
%  \hline
% \end{tabular}
% \end{center}
% \end{table}

% \vspace{-2mm}

% \begin{table}[H]
% \footnotesize
% \caption{\label{tab:LIN error SS param} 500 Perturbed Parameters Linear Approximation Errors without Reduction v.s. with Reduction}
% \begin{center}
% \begin{tabular}{||c|| c c c c c||} 
%  \hline
%   & $\left\lVert u_{lin}^{7}-u_{lin}^{3}\right\rVert$ & $\left\lVert u_{opt}^{7}-u_{lin}^{7}\right\rVert$ & $\left\lVert u_{opt}^{7}-u_{lin}^{3}\right\rVert$ & $\left\lVert u_{opt}^{3}-u_{lin}^{7}\right\rVert$ & $\left\lVert u_{opt}^{3}-u_{lin}^{3}\right\rVert$ \\ [0.5ex] 
%  \hline
%  Mean & 0.0427 & 0.1480 & 0.1447 & 0.1443 & 0.1371 \\ 
%  \hline
%  Median & 0.0437 & 0.1212 & 0.1187 & 0.1260 & 0.1165 \\
%  \hline
% \end{tabular}
% \end{center}
% \end{table}

Table \ref{tab:LIN cost SS param} provides two main results. First, the difference in median values of $J(u_\text{lin}^7)$ and $J(u_\text{lin}^3)$, i.e. \( -9.68 \times 10^{-1} \) and \( -9.42 \times 10^{-1} \) respectively, is minuscule, around a $3\%$ difference. This suggests that parameter reduction is an effective method for reducing computation cost while maintaining accuracy. Secondly, we see that on average a linear approximation does not do a good job of minimizing cost and has many outliers that greatly reduce its accuracy. This is shown by the high mean costs of $2.54\times 10^3$ and $6.83\times 10^2$.

Table \ref{tab:LIN error SS param} provides similar results. We see that the norm difference between approximations using the full parameter space versus the reduced space is extremely small, $4.37 \times 10^{-2}$. This once again suggests that the parameter space reduction does not reduce accuracy of approximations. However, we see that error between any optimal solution and linear approximation is large, such as $\left\lVert u_{opt}^{7}-u_{lin}^{7}\right\rVert$ median of $1.21\times10^{-1}$ which is an 82\% increase from the nominal normed controller. This once again shows that a linear approximation is not effective at approximating an optimal solution.

Tables \ref{tab:IS cost SS param} and \ref{tab:IS error SS param} provide similar data but for our new interpolated step approximation method denoted as $u_{IS}$. The approximations were computed for the same 500 parameter perturbations and for both the full 7 dimensional parameter space and the reduced 3 dimensional parameter space.

% \begin{table}[H]
% \footnotesize
% \caption{\label{tab:IS cost SS param} 500 Perturbed Parameters Interpolated Step Approximation Costs without Reduction v.s. with Reduction}
% \begin{center}
% \begin{tabular}{||c|| c c c||} 
%  \hline
%   & $J\left(u_{opt}^{7}\right)$ & $J\left(u_{IS}^{7}\right)$ & $J\left(u_{IS}^{3}\right)$ \\ [0.5ex] 
%  \hline
%  Mean & 1.3729 & 1.8242 & 1.3374 \\ 
%  \hline
%  Median & -1.6338 & -1.4618 & -1.3891 \\
%  \hline
% \end{tabular}
% \end{center}
% \end{table}

% \begin{table}[H]
% \footnotesize
% \caption{\label{tab:IS error SS param} 500 Perturbed Parameters Interpolated Step Approximation Errors without Reduction v.s. with Reduction}
% \begin{center}
% \begin{tabular}{||c|| c c c c c||} 
%  \hline
%   & $\left\lVert u_{IS}^{7}-u_{IS}^{3}\right\rVert$ & $\left\lVert u_{opt}^{7}-u_{IS}^{7}\right\rVert$ & $\left\lVert u_{opt}^{7}-u_{IS}^{3}\right\rVert$ & $\left\lVert u_{opt}^{3}-u_{IS}^{7}\right\rVert$ & $\left\lVert u_{opt}^{3}-u_{IS}^{3}\right\rVert$ \\ [0.5ex] 
%  \hline
%  Mean & 0.0756 & 0.1103 & 0.1152 & 0.1196 & 0.1038 \\ 
%  \hline
%  Median & 0.0463 & 0.0590 & 0.0542 & 0.0700 & 0.0470 \\
%  \hline
% \end{tabular}
% \end{center}
% \end{table}

\begin{table}[H]
\footnotesize
\caption{\label{tab:IS cost SS param} 500 Perturbed Parameters Interpolated Step Approximation Costs without Reduction v.s. with Reduction}
\begin{center}
\begin{tabular}{||c|| c c c||} 
 \hline
  & $J\left(u_{opt}^{7}\right)$ & $J\left(u_{IS}^{7}\right)$ & $J\left(u_{IS}^{3}\right)$ \\ [0.5ex] 
 \hline
 Mean & $1.37 \times 10^{0}$ & $1.82 \times 10^{0}$ & $1.34 \times 10^{0}$ \\ 
 \hline
 Median & $-1.63 \times 10^{0}$ & $-1.46 \times 10^{0}$ & $-1.39 \times 10^{0}$ \\
 \hline
\end{tabular}
\end{center}
\end{table}

\begin{table}[H]
\footnotesize
\caption{\label{tab:IS error SS param} 500 Perturbed Parameters Interpolated Step Approximation Errors without Reduction v.s. with Reduction}
\begin{center}
\begin{tabular}{||c|| c c c c c||} 
 \hline
  & $\left\lVert u_{IS}^{7}-u_{IS}^{3}\right\rVert$ & $\left\lVert u_{opt}^{7}-u_{IS}^{7}\right\rVert$ & $\left\lVert u_{opt}^{7}-u_{IS}^{3}\right\rVert$ & $\left\lVert u_{opt}^{3}-u_{IS}^{7}\right\rVert$ & $\left\lVert u_{opt}^{3}-u_{IS}^{3}\right\rVert$ \\ [0.5ex] 
 \hline
 Mean & $7.56 \times 10^{-2}$ & $1.10 \times 10^{-1}$ & $1.15 \times 10^{-1}$ & $1.20 \times 10^{-1}$ & $1.04 \times 10^{-1}$ \\ 
 \hline
 Median & $4.63 \times 10^{-2}$ & $5.90 \times 10^{-2}$ & $5.42 \times 10^{-2}$ & $7.00 \times 10^{-2}$ & $4.70 \times 10^{-2}$ \\
 \hline
\end{tabular}
\end{center}
\end{table}

Table \ref{tab:IS cost SS param} shows a similar narrative to the case without interpolation. First, we see that the mean and median are distant because as mentioned earlier, imposing constraints as penalties represents itself as discrepancy. However, this has greatly improved the stability of the result compared to the linear case. The mean and median costs now closely resemble that of the optimized result,  showing that our approximation method does a great job of minimizing cost, just as seen and discussed for figure \ref{fig:ssHist}. Finally, we once again see that the minimal difference between the 7- and 3-parameter cases, median costs of $-1.46$ and $-1.39$ respectively, justifies the use of parameter reduction as it also greatly reduces the computational load of using the interpolation.

Table \ref{tab:IS error SS param} also shows similar norm error results to the linear case. The norm difference between the 7- and 3-parameter approximations, $\left\lVert u_{IS}^{7}-u_{IS}^{3}\right\rVert$ having a median of $4.63 \times 10^{-2}$, is small enough to again imply that parameter reduction is an effective method for reducing computation cost while maintaining accuracy. Additionally, we see much lower errors showing that the interpolated step approximation method accurately replicates the optimal solution while only using precomputed data.

\subsection{Space Shuttle Discussion}
Using the general methods of HDSA, parameter reduction, and our interpolated step method, the shuttle trajectory problem (\ref{eq:ssMIN}) was adjusted accurately for the parameter perturbations without sacrificing robustness. Figure \ref{fig:ssHist} shows right skew that implies the interpolated step approximation minimizes cost and accurately approximates the optimal solution. Table \ref{tab:opt solutions} verifies that parameter reduction does not affect the optimized solutions accuracy. Tables \ref{tab:LIN cost SS param}, \ref{tab:LIN error SS param}, \ref{tab:IS cost SS param}, and \ref{tab:IS error SS param} all show that reducing the parameter space of statistically unimportant parameters does not affect the accuracy of approximated solutions. Additionally, the difference in approximation accuracy from tables \ref{tab:LIN cost SS param} and \ref{tab:LIN error SS param} to tables \ref{tab:IS cost SS param} and \ref{tab:IS error SS param}, verify that our instantaneous interpolated step approximation method does an excellent job at approximating a new optimal controller.

We observe a significant reduction in computation time, as shown in Table \ref{tab:time reduc}. As mentioned earlier, recalculating the optimal solution in-transit is often prohibitive given the lack of any guarantee of a favorable solution in any short period of time. 

\begin{table}[H]
\footnotesize
\caption{\label{tab:time reduc} 500 Perturbed Parameters Approximations Mean Computation Times}
\begin{center}
\begin{tabular}{||c|| c c ||} 
 \hline
  & $\text{time}(u_{opt})$ & $\text{time}(u_{IS})$ \\
 \hline
 Mean (s) & $3.19 \times 10^1$ & $1.12 \times 10^{-1}$ \\
 \hline
\end{tabular}
\end{center}
\end{table}

\noindent Table \ref{tab:time reduc} highlights the considerable decrease in computation time between obtaining the optimal solution and our interpolated step approximation. Note that the primary factor that is limiting the approximation from achieving even faster computation times is the speed of the interpolation function in $\textsc{matlab}$. Additionally, keep in mind that the data in Table \ref{tab:time reduc} was acquired on a local workstation as opposed to the specialized chips and hardware that would be onboard this aircraft problem.

Concisely, our interpolation improves in-situation computation time and approximation stability for both cost and norm error metrics. When coupled with the parameter reduction, the load of precomputation for this method is reduced without a significant sacrifice of accuracy.

\section{Conclusion}

Optimal control under parameter uncertainty restricts the usefulness of the precomputed optimal solution due to in-transit parameter changes. Recalculation of the optimal path is seen to be computationally expensive, time-consuming, and lacks any guarantee of obtaining a favorable solution. The proposed framework of coupling parameter space reduction, an interpolated step approximation, and the use of hyper-differential sensitivity analysis (HDSA), focuses on a quick and accurate in-transit approximation of the optimal controller. For problems in aerospace, by precomputing sensitivity data and solving a simpler problem due to parameter space reduction, instantaneous approximations are able to save the vehicle from inevitable catastrophe. 

As stated in the introduction, we want to emphasize that our method is intended to be used in conjunction with established methods for handling small parameter perturbations such as feedback control. The space shuttle example of Betts verifies the parameter reduction and interpolated-step approximation methods. While this is a limited application, the generality of these methods allows their accuracy to be controlled by the user and allows for more complex application. Stricter methods of applying constraints would also likely limit the instability possible in the approximated solution. It is important to note that this approach is one of many in the toolbox for solving perturbation issues in optimal control problems and thus users should consider other established methods alongside this new approach to determine the most effective and practical solution for the application under consideration. Yet, the proposed approach has applications in aerospace trajectory control, parameter uncertainty calibration, and broad applicability to optimal control problems in which parameter space reduction could improve computation time with little to no reduction in accuracy.  

\section*{Acknowledgments}
As a part of the North Carolina State University's Directed
Research for Undergraduates in Mathematics and
Statistics (DRUMS) program, our research was supported in part by NSF/DMS-2051010, as well as the National Security Agency Grant No. H98230-23-1-0009.

We also thank Alen Alexanderian, Joseph Hart, and Paul Spears for their useful insights and oversight of this project. 

\newpage
\appendix 
\section{Parameters and Variables of the Space Shuttle Problem}\label{appendix:SSvariables}
Given below are the parameters and variables for the 2 DOF space shuttle problem described in section \ref{sec:examples}.

\begin{minipage}{0.5\textwidth}
\begin{table}[H]
\footnotesize
\caption{\label{tab:StaticVar} Static Parameters}
\begin{center}
\begin{adjustbox}{width=0.95\textwidth,center}
\begin{tabular}{||c | c | c||}
    \hline
    Parameter Name & Symbol & Value \\ [0.5ex]
    \hline\hline
    Vehicle Mass & $m$ & $\frac{20300}{32.173}$ ($slug$) \\
    \hline
    Aerodynamic Reference Area & $S$ & $2690$ ($ft^2$) \\
    \hline
    Nominal Atmospheric Density & $\rho_0$ & $0.002378$ ($lb/ft^3$) \\
    \hline
    Scale Height Factor & $h_r$ & $23800$ ($ft$) \\
    \hline
    Radius of the Earth & $Re$ & $20902900$ ($ft$) \\
    \hline
    GM for Earth & $\mu$ & $0.14076539 \times 10^{17}$ ($ft^3/s^2$) \\
    \hline
\end{tabular}
\end{adjustbox}
\end{center}
\end{table} 
\end{minipage}
\begin{minipage}{0.5\textwidth}
\begin{table}[H]
\footnotesize
\caption{\label{tab:DLCoefVar} Dimensionless Drag and Lift Coefficients}
\begin{center}
\begin{adjustbox}{width=0.79\textwidth,center}
\begin{tabular}{||c | c | c||}
    \hline
    Parameter Name & Symbol & Value \\ [0.5ex]
    \hline\hline
    Constant Lift & $a_0$ & $-0.20704$ \\
    \hline
    Linear Lift & $a_1$ & $0.029244$ \\
    \hline
    Constant Drag & $b_0$ & $0.07854$ \\
    \hline
    Linear Drag & $b_1$ & $-0.61592 \times 10^{-2}$ \\
    \hline
    Quadratic Drag & $b_2$ & $0.621408 \times 10^{-3}$ \\
    \hline
\end{tabular}
\end{adjustbox}
\end{center}
\end{table}
\end{minipage}

\begin{table}[H]
\footnotesize
\caption{\label{tab:funcVar} Functional Variables}
\begin{center}
\begin{adjustbox}{width=0.475\textwidth,center}
\begin{tabular}{||c | c | c||}
\hline
Function Name & Variable & Function \\ [0.5ex]
\hline\hline
Atmospheric Density & $\rho(h)$ & $\rho_0\exp{(\frac{h}{h_r})}$ ($lb/ft^3$) \\
\hline
Distance from Earth's Center & $r(h)$ & $Re + h$ ($ft$) \\
\hline
Force of Gravity & $g(h)$ & $\frac{\mu}{r^2}$ ($ft/s^2$) \\
\hline
Angle of Attack & $\hat{u}$ & $u\frac{180}{\pi}$ ($deg$) \\
\hline
Coefficient of Lift & $c_L(\hat{u})$ & $a_0 + a_1\hat{u}$ (-) \\
\hline
Coefficient of Drag & $c_D(\hat{u})$ & $b_0 + b_1\hat{u} + b_2\hat{u}^2$ (-) \\
\hline
Drag Force & $D(c_D,\rho,v)$ & $0.5c_DS\rho v^2$ ($lbft/s^2$) \\
\hline
Lift Force & $L(c_L,\rho,v)$ & $0.5c_LS\rho v^2$ ($lbft/s^2$) \\
\hline
\end{tabular}
\end{adjustbox}
\end{center}
\end{table}
\vspace{1mm}
\section{Space Shuttle Cost Function}\label{appendix:cost}
Given below is the complete cost function for the 2 DOF space shuttle problem. The first term aims to maximize longitudinal distance. The next three penalties enforce terminal conditions and the final four penalties enforce state constraints. The $\beta_i$ are tuning parameters. These tuning parameter values were chosen experimentally and are going to be application dependent. 

% No specific technique was used to choose these in this particular application but in general one could use a generalized L-Curve approach or a hyperparameter tuning technique such as grid search or Bayesian optimzation.

\begin{align*} &J(u) = \frac{-x_{2}(t_{f},u)}{\overline{x_{2}}} \\ & + \beta_1\left( \frac{x_{1}(t_{f},u) - h_{f}}{h_{f}}\right)^{2} +  \beta_2\left( \frac{x_{3}(t_{f},u) - v_{f}}{v_{f}}\right)^{2} + \beta_3\left( \frac{x_{4}(t_{f},u) - g_{f}}{g_{f}}\right)^{2}  \\ & + \beta_4\left(\frac{(\sum_{(x_{1}(t_i) < 0} x_{1}(t_i)^4)}{\overline{x_{1}}^4}\right) + \beta_5\left(\frac{( \sum_{(x_{3}(t_i) < 1} (x_{3}(t_i)-1)^4))}{\overline{x_{3}}^4}\right) \\ &+ \beta_6\left(\frac{(\sum_{(x_{4}(t_i) < \frac{-89 \pi}{180}} (x_{4}(t_i)+\frac{89 \pi}{180})^4))}{\overline{x_{4}}^4}\right)   + \beta_7\left(\frac{(\sum_{x_{4}(t_{i})> \frac{89 \pi}{180}} (x_{4}(t_i)-\frac{89 \pi}{180})^4)}{\overline{x_{4}}^4}\right)\end{align*}

\section{Sensitivity Equations}\label{appendix:senseq}
The sensitivity equations compute $\frac{dJ}{du_{i}}$ using a manipulation of chain rule that allows the exploitation of ODE solvers to make these derivatives more accurate than finite difference or complex step when calculated numerically. From the defining dynamical systems, we have $\frac{dx}{dt} = F$. Differentiating this gives 
\begin{equation*}
\frac{\partial^{2}x}{\partial u_{i}\partial t} = \frac{\partial F}{\partial x} \frac{\partial x}{\partial u_{i}} + \frac{\partial F}{\partial u} \frac{\partial u}{\partial u_{i}}
\end{equation*}

Since the mixed partials are continuous, we may apply Clairaut's theorem that $\frac{\partial^{2}x}{\partial u_{i}\partial t} = \frac{\partial ^{2}x}{\partial t\partial u_{i}}$, establishing an ODE against time from which we solve for $\frac{\partial x}{\partial u_{i}}$. In the context of the space shuttle, we can set up this ODE as 
\begin{equation*} 
\frac{dF}{dx} = \begin{bmatrix}
 0&0& \sin(x_4)& x_{3}\cos(x_4)\\
 \frac{-x_{3} \cos(x_4)}{r^{2}}&0&\frac{\cos(x_4)}{r}& -\frac{x_{3}}{r}\sin(x_4)\\ 
 \frac{-D}{h_r} + \frac{2g\sin(x_4)}{r}&0&\frac{2D}{x_3}& -g\cos(x_4)\\
 \frac{-L}{x_{3}mh_{r}} + cos(x_4)(\frac{-x_3}{r^{2}} + \frac{2g}{x_3 r})&0&\frac{L}{mx_{3}^{2}} + \cos(x_{4})(\frac{1}{r} + \frac{g}{x_{3}^{2}})& \sin(x_4)(\frac{x_{3}}{r} - \frac{g}{x_{3}})
\end{bmatrix}
\end{equation*}
and 
\begin{equation*}
\frac{dF}{du} = \begin{bmatrix} 
0 \\ 
0 \\  
\frac{-1}{2m}(S\rho_{0}e^{-x_{1}/h_{r}}x_{3}^{2})(b_1+2b_2 \hat{u})(\frac{180}{\pi})\\
\frac{1}{m}(S\rho_{0}e^{-x_1/h_{r}}x_{3}^{2})(a_1)(\frac{180}{\pi})
\end{bmatrix}
\end{equation*}

% Next recall the cost function from appendix \ref{appendix:cost}. We differentiate this with respect to $u_{i}$. First, consider $-x_{2}(t_{f},u)$. By the Lebesgue differentiation theorem we may write $$-x_{2}(t_{f},u) = \lim_{n \to \infty} n\int_{t_{f} - \frac{1}{n}}^{t_{f}} -x_{2}(t,u)dt$$ Along with the continuity of $x_2$ and its partials with respect to $u_i$ and $t$, this representation gives by the differentiable limit theorem that 
% \begin{equation*}
% \frac{\partial}{\partial u_{i}} (-x_{2}(t_{f},u)) = \lim_{n \to \infty} n \frac{\partial}{\partial u_{i}}\int_{t_{f} - \frac{1}{n}}^{t_{f}} -x_{2}(t,u)dt
% \end{equation*}

% The Leibniz integration rule allows for differentiation under the integral sign, so that \begin{equation*}
% \frac{\partial}{\partial u_{i}} (-x_{2}(t_{f},u)) = \lim_{n \to \infty} n \int_{t_{f} - \frac{1}{n}}^{t_{f}} -\frac{\partial x_{2}(t,u)}{\partial u_{i}} dt = -\frac{dx_{2}}{du_{i}}(t_f,u) 
% \end{equation*}

% This computation works similarly for the means of the $x_{i}$. We may also note that.
Through applications of the Lebesgue differentiation theorem and Leibniz integration rule, as well as considering the piecewise portions of the cost, we arrive at
\begin{align*} &\frac{\partial J}{\partial u_i} =  \frac{\overline{x_{2}}(-\frac{dx_{2}}{du_{i}}(t_{f},u)) +x_{2}(t_{f},u)\overline{\frac{\partial dx_{2} }{\partial u_{i} }}) }{(\overline{x_{2}})^{2}} \\
&+2\beta_{1}\frac{(x_{1}(t_{f},u) - h_{f})}{h_{f}^{2}} \frac{\partial x_{1}}{\partial u_{i}}(t_{f},u)  +2\beta_{2}\frac{(x_{3}(t_{f},u) - v_{f})}{v_{f}^{2}} \frac{\partial x_{3}}{\partial u_{i}}(t_{f},u) + 2\beta_{3}\frac{(x_{4}(t_{f},u) - g_{f})}{g_{f}^{2}} \frac{\partial x_{4}}{\partial u_{i}}(t_{f},u) \\ 
&+ 4\beta_4\left( \frac{(\overline{x_1}^{4}(\sum_{x_{1}(t_{i})<0}4x_{1}(t_{i})^3\frac{\partial x_{1}}{\partial u_i}(t_{i},u)) -  (\sum_{x_{1}(t_{i})<0}x_{1}(t_i)^{4})(4\overline{x_1}^3 \overline{\frac{\partial x_1}{\partial u_i}})}{\overline{x_{1}}^{8}}\right) \\ 
&+ 4\beta_5 \left( \frac{(\overline{x_3}^{4}(\sum_{x_{3}(t_{i})<1}4(x_{3}(t_{i})-1)^3\frac{\partial x_{3}}{\partial u_i}(t_{i},u)) -  (\sum_{x_{3}(t_{i})<1}(x_{3}(t_i)-1)^{4})(4\overline{x_3}^3 \overline{\frac{\partial x_3}{\partial u_i}})}{\overline{x_{3}}^{8}}\right)\\ 
&+ 4\beta_6 \left( \frac{(\overline{x_4}^{4}(\sum_{x_{4}(t_{i})<\frac{-89 \pi}{180}}4(x_{4}(t_{i})+\frac{89 \pi}{180})^3\frac{\partial x_{4}}{\partial u_i}(t_{i},u)) -  (\sum_{x_{4}(t_{i})<\frac{-89 \pi}{180}}(x_{4}(t_i)+\frac{89 \pi}{180})^{4})(4\overline{x_4}^3 \overline{\frac{\partial x_4}{\partial u_i}})}{\overline{x_{4}}^{8}}\right)\\ 
&+ 4\beta_7 \left( \frac{(\overline{x_4}^{4}(\sum_{x_{4}(t_{i})>\frac{89 \pi}{180}}4(x_{4}(t_{i})-\frac{89\pi}{180})^3\frac{\partial x_{4}}{\partial u_i}(t_{i},u)) -  (\sum_{x_{4}(t_{i})>\frac{89\pi}{180}}(x_{4}(t_i)-\frac{89\pi}{180})^{4})(4\overline{x_4}^3 \overline{\frac{\partial x_4}{\partial u_i}})}{\overline{x_{4}}^{8}}\right) \end{align*}

\pagebreak
\bibliographystyle{siamplain}
\bibliography{references}

\end{document}